\input amstex
  \documentstyle{amsppt}
  \NoBlackBoxes
  \magnification=1200
  \NoRunningHeads
  \redefine\c{\Bbb C}
  \topmatter
  \title Nonuniformizable skew cylinders. A counterexample to the simultaneous
  uniformization problem\endtitle
  \author A.A.Glutsyuk\endauthor
  \address \ \ Independent University \ of Moscow, \ \ Bolshoi Vlasievskii Pereulok 11,
  121002 Moscow Russia\endaddress
  \address Steklov Mathematical Institute, Moscow\endaddress

  \address Present Address: Institut des Hautes \'Etudes Scientifiques,
  Le Bois-Marie, 35 Route de Chartres, 91440 Bures-sur-Yvette, France\endaddress

  \thanks Research supported by CRDF grant RM1-229, by
  INTAS grant 93-0570-ext,
  by Russian Foundation for Basic Research (RFBR) grant 98-01-00455, by State
  Scientific Fellowship of Russian Academy of Sciences for young scientists,
  by the European Post-Doctoral Institut joint fellowship of Max-Planck Institut f\"ur
  Mathematik (Bonn) and IH\'ES (Bures-sur-Yvette, France)\endthanks
  
  \abstract At the end of 1960$^{ths}$ Yu.S.Ilyashenko
  stated the
  problem: is it true that for any one-dimensional holomorphic foliation
  with singularities on a Stein manifold leaves intersecting
  a transversal disc can be uniformized by family of simply connected domains
  in $\overline{\Bbb C}$
  so that the uniformization function would depend holomorphically
  on the transversal parameter?
  In the present paper we give a two-dimensional counterexample.
  This together with the previous result of Ilyashenko (Lemma 1)
  implies existence of a counterexample given by a foliation
  on affine (projective) algebraic surface by level curves of a rational function
  with singularities
  deleted. This implies that {\it Bers' simultaneous uniformization theorem
  \cite{9} for topologically trivial holomorphic fibrations by compact Riemann
  surfaces does not extend to general fibrations by compact Riemann surfaces with
  singularities.}
  \endabstract
  \endtopmatter
  \document
  
  \subhead 1. Skew cylinders and simultaneous uniformization\endsubhead
  
  By $T_g$, $g\geq2$, denote the Teichm\"uller space of closed Riemann surfaces
  of genus $g$. The tautological fibration over $T_g$ is the union of all
  the Riemann surfaces: the fiber over each point of $T_g$ is the corresponding
  Riemann  surface.
  The tautological fibration admits a natural complex structure, and so does
  the union of the universal coverings of its fibers. The latter is fibered over
  $T_g$ with fibers conformally equivalent to disc.
  
  In 1960 L.Bers \cite{9} proved the classical simultaneous uniformization
  theorem saying that fibers of the tautological fibration
  admit a uniformization that depends holomorphically on the Teichm\"uller space
  parameter. Namely, he proved that the union of the universal coverings of
  the fibers is biholomorphically fiberwise
  equivalent to a subdomain in the direct product $\overline{\Bbb C}\times T_g$.
  This theorem implies the analogous statement for any topologically
  trivial rational fibration by algebraic curves. In particular, for any
  two-dimensional smooth projective surface
  $S$ fibered by algebraic curves with singularities and any transversal disc
  $D\subset S$
  intersecting only nonsingular fibers the union of the universal coverings of
  the fibers intersecting $D$ is  biholomorphically fiberwise equivalent to a
  subdomain in the direct product $\overline{\Bbb C}\times D$.
  
  It appears that the generalization of the last statement to topologically
  nontrivial fibrations by algebraic curves {\it with singularities} is wrong
  in general.
  
   \proclaim{Theorem 1} There exists a smooth algebraic surface $S$ (of dimension 2)
   in $\Bbb C^n$ ($\Bbb P^n$) fibered by algebraic curves with singularities
   such that the family of Riemann surfaces thus obtained does not admit
   simultaneous uniformization by family of simply connected domains in
   $\overline{\Bbb C}$ such that the uniformizing function depends
   holomorphically on the parameter. (The number $n$ may be chosen to be equal
   to 5.)
  
   More precisely, there exists
  a polynomial $P$ in $\c^n$, $P|_S\not\equiv const$, and an inverse
  mapping $i:D\to S$ of the unit disc $D$ in the image of $P$, $P\circ i=Id$,
  with the following properties:
  
  \item{1)} The image $i(D)$ is transversal to the foliation $P=const$ by level
  curves of $P$.
  
  \item{2)} Let $M$ be the union of the universal coverings of the leaves
  intersecting $i(D)$ of the foliation $P=const$ equipped with the natural
  complex structure. There is no biholomorphic mapping of $M$ onto a domain
  in $\overline{\Bbb C}\times D$ that forms a commutative diagram with $P$
  and the standard direct product projection.
  \endproclaim
  
  Theorem 1 is proved at the end of the Section.
  
  \definition{Definition 1} Let $D$ be a simply-connected domain in complex line,
   $M$ be a two-dimensional complex
  manifold, $p:M\to D$ be a proper holomorphic surjection having nonzero
  derivative. We say that
  the triple $(M,p,D)$ is a {\it skew cylinder} with the base $D$ and the total
  space $M$, if
  
  1) all the level sets of the map $p$ are simply connected holomorphic curves;
  
  2) the "fibration" $M$ admits a section, i.e., a
  holomorphic mapping $i:D\to M$ such that $p\circ i=Id$.
  \enddefinition
  
  \example{Example 1} Let $S$ be a smooth affine (projective) two-dimensional
  complex surface, $P$ be a nonconstant rational function, $D$ be a disc in the
  image of $P$,
  $i:D\to S$ be an inverse mapping transversal to the foliation $P=const$,
  $P\circ i=Id$. Let $M$ be the union of the
  universal coverings of the leaves intersecting $D$ equipped with the natural
  complex structure. Then the manifold $M$ equipped with the standard
  projection $P:M\to D$ is a skew cylinder. A skew cylinder thus constructed is
  called {\it affine (projective) algebraic}.
  \endexample
  
  \definition{Definition 2} A skew cylinder $(M,p,D)$ is said to be
  {\it uniformizable}, if there exists a biholomorphic mapping
  $u: M\to \overline{\Bbb C}\times D$ (not necessarily "onto") that forms a
  commutative diagram with the projections to $D$.
  \enddefinition
  
  \definition{Definition 3} A skew cylinder is said to be {\it Stein}, if its
  total space is Stein.
  \enddefinition
  
  The main result of the paper is the following
  
  \proclaim{Theorem 2} There exists a nonuniformizable Stein skew cylinder.
  \endproclaim
  
  Theorem 2 is proved in Subsection 2.4.
  
  \definition{Definition 4} Let $(M,p,D)$ be a skew cylinder, $B\subset M$ be
  its
  subdomain. Then $B$ is called a {\it subcylinder}, if the triple $(B,p,p(B))$ is
  a skew cylinder.
  \enddefinition
  
  \definition{Definition 5} Two skew cylinders are said to be {\it equivalent}, if
  there exists a biholomorphism of the total space of the one onto that of the other
  that forms a commutative diagram with the projections.
  \enddefinition
  
  \proclaim{Lemma 1 (Yu.S.Ilyashenko)}\footnote{From unpublished paper by
  Yu.S.Ilyashenko (late 1960$^{ths}$)} Any compact subcylinder of a Stein skew
  cylinder is equivalent to a subcylinder of an affine (projective) algebraic
  skew cylinder.
  The dimension of the ambient affine (projective) space of the corresponding
  surface (see Example 1) may  be chosen to be equal to 5.
  \endproclaim
  
  Lemma 1 is proved in Subsection 2.5.
  
  \proclaim{Proposition 1} Let a Stein skew cylinder be exhausted by increasing
  sequence of uniformizable subcylinders. Then it is uniformizable.
  \endproclaim
  
  A version of Proposition 1 was proved by Yu.S.Ilyashenko in \cite{6}.
  The proof of its present version is a modification of Ilyashenko's proof
  (his argument that uses normality of the space of normalized $\Bbb C$- valued
  univalent functions in unit disc \cite{6} is replaced by analogous one using normality of
  the space of normalized $\overline{\Bbb C}$- valued univalent functions $\phi$:
  $\phi(0)=\phi''(0)=0$, $\phi'(0)=1$).
  
  \demo{Proof of Theorem 1} By Theorem 2, there exists a
  nonuniformizable Stein skew cylinder. By Proposition 1, it has a nonuniformizable compact
  subcylinder. By Lemma 1, this subcylinder is equivalent to a
  subcylinder of an affine (projective) algebraic skew cylinder. The last
  subcylinder is nonuniformizable as well. The same holds for the ambient
  algebraic cylinder. This proves Theorem 1.
  \enddemo
  
  \head 2. Nonuniformizability\endhead
  
  In this Section we prove Theorem 2 (Subsection 2.4) and Lemma 1 (Subsection 2.5).
  
  In Subsections (2.1-2.3) we give a survey of previous results and state an
  open question (Subsection 2.3).
  
  \subhead 2.1. Previous results for algebraic skew cylinders\endsubhead
  
  In 1973 Yu.S.Ilyashenko \cite{5} extended Bers'  simultaneous uniformization
  theorem to the tautological fibration with a unique singular fiber
   over disc embedded in the Deligne-Mumford compactification
  of the moduli space of Riemann surfaces.
  This implies the positive solution of the simultaneous uniformization problem
  for foliation with singularities
  by level curves of rational function and leaves over a disc in its image
  under the assumption that
  {\it only one of these leaves is singular} and has at most simple double
  point singularities.
  
  \subhead 2.2. Previous results for general skew cylinders\endsubhead
  
  It is easy to construct nonstein nonuniformizable skew cylinders.
  
  \example{Example 2} Let $D$ be unit disc, $f:D\to\c$ be arbitrary 
  nonholomorphic
  function (say, $f(z)=\bar z$). Let $M$ be the universal covering over the
  complement to the graph of $f$
  in the direct product $\c\times D$. The manifold $M$ admits a natural
  structure of skew cylinder, and the latter is nonuniformizable.
  \endexample
  
  Other examples  of nonstein nonuniformizable skew cylinders with fibers
  conformally equivalent to $\c$ may be found in \cite{6}.
  
  Yu.S.Ilyashenko \cite{6} proposed a conjecture that any Stein skew cylinder is
  uniformizable.
  In 1969 T.Nishino \cite{10} independently proved the positive answer for skew
  cylinders with fibers $\c$.
  
   A.A.Shcherbakov \cite{14} proved that a Stein skew cylinder can be exhausted
   by a growing sequence of compact subcylinders with smooth
   strictly pseudoconvex boundaries.
  
  Proposition 1 together with Theorem 2 imply the following
  
  \proclaim{Corollary 1} There exists a compact skew cylinder with a smooth
  strictly pseudoconvex boundary that is nonuniformizable.
  \endproclaim
  
  \subhead 2.3. Foliations on $\Bbb P^2$ and uniformization problem for them
  \endsubhead
  
  It is well known that in the generic case the arrangement of phase
  curves of a polynomial vector field in $\c^n$ (or more generally, leaves of
  foliation with singularities by analytic curves) is quite complicated.
  For example, there is an open domain $U$ in the space of polynomial vector fields
  in $\c^2$ such that each phase curve of a generic polynomial vector field
  from $U$ is everywhere dense in the phase space \cite{1,2,3}. A statement of this
  type follows from a theorem of I.Nakai \cite{16}. Recently
  F.Loray and J.Rebelo proved in their joint paper \cite{13} that there is an
  open subset $V$ in the space
  of one-dimensional holomorphic foliations with singularities on $\Bbb P^n$
  of fixed degree greater than 1 such that for any foliation from $V$ each leaf
  is dense.
  
  
   \definition{Definition 6} Let $S$ be a two-dimensional complex manifold,
   $F$ be a one-dimensional
  holomorphic foliation with isolated singularities on $S$. Let $D$ be an
  embedded disc in $S$ transversal to $F$. The {\it universal covering manifold}
  associated to the triple $(S,F,D)$ is the union of the universal coverings
  over the leaves intersecting $D$
   with marked points in $D$. More precisely, it is the set of all the triples
   consisting of a point $z\in D$, a point $z'$ of the leaf containing $z$,
   and a homotopy class of a path  connecting $z$ to $z'$ in the leaf.
   \enddefinition
  
   \remark{Remark 1} In the condition of the previous Definition let $S$ be Stein
   (or projective).
   Then for any transversal disc $D$ the corresponding universal covering manifold
   admits a natural structure of complex manifold; it is Stein if $S$ is Stein.
   These statements hold in any dimension for one-dimensional holomorphic
   foliations with singular sets of complex codimension at least 2. For Stein $S$
    they were proved by Ilyashenko (\cite{4}, \cite{15}). For projective $S$ the first statement
     was proved by myself and a little later by E.M.Chirka (in unpublished
     papers). In both cases the universal covering manifold has a natural skew
     cylinder structure with the base $D$.
   \endremark
  
  \proclaim{Question 1} Does Theorem 1 hold with $S=\c^2$ or $\Bbb P^2$?
  Does there exist a polynomial (holomorphic) vector field
  in $\c^2$ (or a one-dimensional holomorphic foliation with isolated
  singularities on $\Bbb{P}^2$) and a transversal disc such that the
  corresponding universal covering manifold is not uniformizable?
  \endproclaim
  
  \subhead 2.4. Proof of Theorem 2\endsubhead
  
  Let $D$ be unit disc in complex line.
  For the proof of Theorem 2 we show that there exists a closed subset
  $K\subset\Bbb C\times D$ fibered over $D$ by discs (an infinite number
  of them degenerates into single points)
  such that the universal covering over its complement $\Bbb C\times D\setminus K$
  (equipped with the canonic projection to $D$) is a nonuniformizable Stein skew
  cylinder. We prove this statement for the following set $K$
  constructed by Bo Berndtsson and T.J.Ransford in their joint paper \cite{7}.
  
  \proclaim{Theorem 3 \cite{7}} Let $D$ be unit disc in complex line
  with the coordinate $z$,
  $p$ be the standard projection to $D$ of the direct product $\Bbb C\times D$.
  Let $E_+=\{\frac12,\frac{n}{2n+1}\}_{n\in\Bbb N}, \ E_-=-E_+\subset D$.
  There exists a closed subset $K\subset \Bbb C\times\overline D$ such that
  
  1) the complement $M'=\Bbb C\times D\setminus K$ is pseudoconvex;
  
  2) for any $z\not\in E_+\cup E_-$
  the fiber $K\cap p^{-1}(z)$ is a disc;
  
  3) for any $z\in E_+$ $K\cap p^{-1}(z)=0\times z$;
  
  4) for any $z\in E_-$ $K\cap p^{-1}(z)=1\times z$.
  \endproclaim
  
  For the completeness of presentation, we recall the construction of the
  set $K$ from \cite{7}. Let $w$ be the coordinate in the
  fiber $\Bbb C$ in the direct product $\Bbb C\times D$.
  Let $u(z)=\ln|z-\frac12|+\ln|z+\frac12|+\sum_{n=1}^{+\infty}2^{-n}(\ln|z-\frac n{2n+1}|+
  \ln|z+\frac n{2n+1}|)$, $A\in\Bbb R^+$. The function $u$
  is harmonic and is equal to $-\infty$ at $E_{\pm}$.
  Let $\psi:D\to\Bbb C$
  be a $C^{\infty}$ function with bounded derivatives (up to the second order)
  that is constant in a neighborhood of each one of the sets $E_{\pm}$ so that
  $\psi|_{E_+}=0$, $\psi|_{E_-}=1$.
  Define
  $$K=\{|w-\psi(z)|\leq e^{u(z)+|z|^2+A}\}.$$
  The fibers of $K$ over $E_+$ ($E_-$) are single points where
  the coordinate $w$ is equal to 0 and 1 respectively.
  If $A$ is large enough, then $\Bbb C\times D\setminus K$
  is pseudoconvex. The proof of this statement (presented in \cite{7}) is a
  straightforward calculation of the Levi form together
   with the argument on approximation of the harmonic
  function $u$ by decreasing sequence of smooth subharmonic functions.
  
  The nonuniformizable skew cylinder we are looking for
  is the universal covering over the complement $M'=\Bbb C\times D\setminus K$
  (denote this covering by $M$). Equivalently, it is the universal covering
  manifold for the foliation on $M'$ by the level curves of the projection
  to $D$ (all these level curves intersect one and the same transversal disc
  $N\times D$ corresponding to $N$ large enough). Indeed, $M$ is Stein.
  This follows from pseudoconvexity of the complement
  $\Bbb C\times D\setminus K$ and the classical theorem due to Stein \cite{12}
  saying that a covering over a Stein manifold is Stein.
   Let us prove that $M$ is nonuniformizable by contradiction. Suppose the
   contrary, i.e., there exists a
  uniformization $u: M\to \overline{\Bbb C}\times D$. Let $\pi:\overline\c\times
  D\to\overline\c$ be the standard projection, $f=\pi\circ u$ be
  the corresponding coordinate component of  $u$.
   The fibers of the cylinder $M$ over $E_{\pm}$ are conformally
  equivalent to complex plane. Let $w$ be the coordinate
  in $\Bbb C$ (we consider it as a coordinate on
  $M'\subset\Bbb C\times D$). Consider the chart
  $\ln w$ on the fibers of $M$. It is a well-defined
  1-to-1 chart on the fibers over $E_+$, and this is not the case for
  the fibers over $E_-$, where this chart is multivalued and
  has branch points. Therefore, for any $z\in E_+$ the restriction to the fiber
  over $z$ of the function $f$  (which is univalent) is M\"obius in the chart
  $\ln w$, and this is not the case for the fibers over $E_-$.
  
  Let $Sf$ be the Schwartzian
  derivative of $f$ in the coordinate $\ln w$. It is a holomorphic
  function on $M$. It vanishes identically on all the fibers of $M$ over $E_+$ and
  on no fiber over no point in $E_-$ by the previous statement.
  The first one of the two last
  statements implies that $Sf\equiv0$ (the set $E_+$ contains the limit point
  $\frac12$), which contradicts the second statement. This proves Theorem 2.
  
  \subhead 2.5. Proof of Lemma 1\endsubhead
  \redefine\t#1{\widetilde#1}
  The proof of Lemma 1 presented below essentially repeats that due to Yu.S.Ilyashenko
  \footnote{From unpublished paper by Yu.S.Ilyashenko (late 1960$^{ths}$)}.
  
  Let $(M,p,D)$ be a Stein skew cylinder.
   The manifold $M$ is Stein. We consider that it is embedded as a submanifold
   in $\Bbb C^n$. (Without loss of generality one can consider that $n=5$
  by Bishop-Narasimhan embedding theorem.) Let $B\subset M$ be a
  subcylinder with compact closure.
  Without loss of generality we consider that its base $p(B)$ is unit disc.
  We will construct a smooth affine algebraic surface $S\subset\c^n$ and a
    polynomial $P$ such that there exists a biholomorphism $h:B\to S$ (not onto)
    that forms
    a commutative diagram with $p$ and $P$. Consider the foliation $P=const$
    (with singularities) by level curves of $P$, or its extension up to a
    foliation with isolated singularities on the projective closure of $S$ 
    (the surface $S$ constructed below may be chosen to have a smooth projective 
    closure). Let $M_B$ be the union of the
    universal coverings of the leaves intersecting the image $h(B)$. Then
    $M_B$ is an algebraic skew cylinder from Lemma 1 we are looking for.

  The surface $M$ is a complete intersection: there exist $n-2$ holomorphic
  functions $f_1,\dots,f_{n-2}$ in $\c^n$ such that $M$ is the transversal
  intersection of their zero level hypersurfaces
  (corollary 1.5 from the paper \cite{8}). Let $i:D\to M$ be the section of $M$.
  The projection $p:M\to D$ extends up to a holomorphic function
(still denoted by $p$) on $\c^n$ by classical theorem on extension of 
analytic function on submanifold of Stein manifold \cite{11}.
  
    To construct the surface $S$ and the polynomial $P$, we use the fact that
    the functions $f_i$, $i=1,\dots,n-2$, and $p$  can be approximated by
    polynomials uniformly
  on each ball in $\Bbb C^n$. Let us fix a ball $\tilde B$ containing $\overline B$ and
  denote the corresponding approximating polynomials by $F_i$ and $P$ respectively.
  Let us denote by $S$ the irreducible component of the conjoint zero
  set of the polynomials $F_i$ that approaches $B$, as $F_i\to f_i$, $P\to p$.
  Let us prove that there exists a
  biholomorphic mapping $h:\overline B\to S$ that forms a commutative diagram with
  the projections $p$ and $P$. To construct this biholomorphism $h$,
  we use the following
  
  \proclaim{Proposition 2} Let $(M,p,D)$ be a Stein skew cylinder, $i:D\to M$ be
  its section (see Definition 1). There exists a holomorphic
  function $f$ on $M$ with nonvanishing derivative along each fiber $p^{-1}(z)$,
  $z\in D$, such that $f|_{i(D)}=0$.
  \endproclaim
  
  \demo{Proof} Let $(M,p,D)$, $i$ be as in Proposition 2. Let us construct a
  function $f$ on $M$ that satisfies the statements Proposition 2.
  
  A generic holomorphic function on $M$ has nonidentically-vanishing derivative
  along each fiber (since $M$ is Stein). Let us fix such a function $g(x)$ that in
  addition vanishes at the image $i(D)$ of the section (one can achieve this by
  changing $g$ to the function $g-g\circ i\circ p$).
  By $\Sigma$ denote the subset of points in $M$ where the restriction of the
  differential $dg$ vanishes on the line tangent to the fiber. In the case,
  when $\Sigma$ is empty, $g$ is a function we are looking for.
  Suppose $\Sigma\neq\emptyset$. Then $\Sigma$ is a hypersurface such that no its
  irreducible component is contained in a fiber. By
  definition, the {\it multiplicity} of its irreducible component $\Sigma'$ is
  said to be the order of zero of the restriction to the fiber of the
  differential $dg$ at a generic point
  of $\Sigma'$. There exists a holomorphic function
  on $M$ that vanishes at $\Sigma$ with the prescribed multiplicities at
  its irreducible components (let us fix this function and denote it by $F$).
  This follows from the theorem saying that any holomorphic line bundle over a
  contractible Stein manifold is trivial. (This is a corollary of the
  classical theorem on triviality of the first
  $\overline{\partial}$- cohomology group for Stein maniflods \cite{11}.)
  By construction, the ratio $\frac{dg}{F}$ is a nowhere vanishing 1-form
  holomorphic on each fiber. The function
  $$f(x)=\int_{i\circ p(x)}^x\frac{dg}{F}$$
  on $M$ is a one we are looking for (the integration is made along a path in
  the fiber containing $x$ that starts at $i\circ p(x)$ and ends at $x$).
  Proposition 2 is proved.\enddemo
  
  Let us construct the fiberwise biholomorphism $h:\overline B\to S$. To do this,
  consider a function $f$ from Proposition 2 and its extension to
  $\Bbb C^n$ (denoted by the same symbol $f$). The function $f$ is locally (but
  not globally) univalent on each fiber. By construction, the level curves
  of $f$ in $\Bbb C^n$ are transversal to the fibers of $M$. The biholomorphism $h$
  we are looking for is defined by the following pair of equations:
  $p(x)=P(h(x))$, $f(x)=f(h(x))$. It is well defined provided that
  the above approximations of the surface $M\cap\tilde B$ and the function
  $p|_{\tilde B}$ by $S\cap\tilde B$  and $P$ respectively are accurate enough
  and the more accurate they are, the more close $h$ is
  to identity. It forms a commutative diagram with $p$ and $P$ by
  construction. Lemma 1 is proved.
  
  \head Acknowledgements\endhead
  
  I am grateful to Yu.S.Ilyashenko, who attracted my attention to the problem.
  I am grateful to G.M.Henkin, who informed me about
  the paper \cite{7} (one of results from \cite{7} became the base for the proof).
  I wish to thank both them and M.Gromov, R.Gunning, C.D.Hill, N.G.Kruzhilin,
  I.Lieb, S.Nemirovski, A.A.Shcherbakov for helpful discussions.
   The paper was written
  when I was visiting the Institut des Hautes \'Etudes Scientifiques
  (Bures-sur-Yvette, France). I wish to thank the Institute for hospitality
  and support.
  
  \head References\endhead
  
   1. Arnold, V.I. and Ilyashenko Yu.S., Ordinary Differential Equations. -
   Itogi Nauki i Techniki, Contemporary Problems in Mathematics, Fundamental
   Directions, Moscow: VINITI, 1985, vol.1, part 1.
  
   2. Muller, B., On the Density of the solutions of Certain Differential
   Equations in $\Bbb{CP}^n$. - Mat. Sb., 1975, vol.98, no. 3, pp. 363-377.
  
   3. Khudai-Verenov, M.O., On a property of the Solutions of a Differential
   Equation. - Mat. Sb., 1962, vol.56 (98), no.3, p.301-308.
  
   4. Ilyashenko, Yu.S., Foliations by analytic curves. -
  Mat. Sb., {\bf88}, No.~4, 558--577 (1972).
  
   5. Ilyashenko Yu.S., Nondegenerate $B$- groups. - Doklady Acad. Nauk SSSR
  208 (1973), 1020-1022. (English translation in Soviet Math. Dokl. 14 (1973),
  207-210.)
  
   6. Ilyashenko, Yu.S. and Shcherbakov, A.A., Skew Cylinders and Simultaneous
   Uniformization. - Proc. of Steklov Math. Inst., 1996, vol.213.
  
   7. Berndtsson, Bo; Ransford, T.J. Analytic multifunctions, the
   $\overline\partial$- equation, and a proof of the corona theorem.- Pacific J.Math. 124
  (1986) No 1, 57-72.
  
    8. Forster, O.;  B\u anic\u a, C. \ \ Complete \ intersections \ \ \ in\  Stein\ Manifolds. -
    Manuscripta Math. 37 (1982), no. 3, 343-356.
  
    9. Bers, L.; Simultaneous uniformization. - Bull. Amer. Math. Society.
    66 (1960) 94-97.
  
    10. Nishino, T.; Nouvelles recherches sur les fonctions enti\`eres
    de plusieurs variables complexes (II). Fonctions enti\`eres qui se
    reduisent \`a celles d'une variable. - J. Math. Kyoto Univ. 9-2 (1969),
    221-274.
  
    11. Gunning, R. Introduction to holomorphic functions of several complex
    variables. Vol.III. Homological theory. - The Wadsworth $\&$ Brooks$\slash$Cole
    Mathematics Series.  Wadsworth $\&$ Brooks$\slash$Cole Advanced Books $\&$
    Software, Monterey, CA, 1990.
  
    12. Stein, K. \"Uberlagerungen holomorph-vollst\"andiger komplexer R\"aume. -
    Arch. Math. 1956, 7 No 5, 354-361.
  
    13. Loray, F.; Rebelo, J. Stably chaotic rational vector fields on $\Bbb{CP}^n$. -
    To appear.
  
    14. Shcherbakov, A.A. The exhaustion method for skew cylinders, - to appear in
    Algebra i Analiz.
  
    15. Ilyashenko, Yu.S. Covering manifolds for analytic families of leaves of
    foliations by analytic curves. - Topol. Methods Nonlinear Anal., 11 (1998),
    no.2, 361-373.
  
    16. Nakai, I. Separatrices for nonsolvable dynamics on $\Bbb C,0$. -
    Ann. Inst. Fourier, Grenoble 44, 2 (1994), pp. 569-599.
  

  \enddocument